\magnification\magstep 0

\documentstyle{amsppt}
\def\a{\alpha} \def\b{\beta} \def\d{\delta}
\def\g{\gamma}  \def\k{\kappa}
\def\l{\lambda} \def\s{\sigma}
 \def\sm{\setminus} \def\o{\omega}
\def\w{\omega}
\def\t{\tilde}
\def\to{\rightarrow}  \def\imp{\Rightarrow}
\def\r{\upharpoonright}  \def\cl{\overline}
  \def\0{\emptyset}
 
\def\<{\langle} \def\>{\rangle}
\topmatter \nologo
\title A note on D-spaces
\endtitle

\author Gary Gruenhage \endauthor
\address Dept. of Mathematics, Auburn University, Auburn, AL 36849
\endaddress
\email garyg\@auburn.edu \endemail

\subjclass 54D20 \endsubjclass
\thanks I would like to dedicate this paper to my colleague and
friend A.V. Arhangel'skii on the occasion of his $65^{th}$
birthday. \endthanks
\thanks Research partially supported by
National Science Foundation grant DMS-0072269 \endthanks

\abstract  We introduce notions of nearly good relations and
$N$-sticky modulo a relation as tools for proving  that spaces are
$D$-spaces. As a corollary to general results about such
relations, we show that $C_p(X)$ is hereditarily a $D$-space
whenever $X$ is a Lindel\"of $\Sigma$-space. This answers a
question of Matveev, and improves a result of Buzyakova, who
proved the same result for $X$ compact.

We also prove that if a space $X$ is the union of finitely many
D-spaces, and has countable extent, then X is linearly Lindel\"of.
It follows that if X is in addition countably compact, then X must
be compact.  We also show that Corson compact spaces are
hereditarily $D$-spaces.  These last two results answer recent
questions of Arhangel'skii.  Finally, we answer a question of van
Douwen by showing that a perfectly normal collectionwise-normal
non-paracompact space constructed by R. Pol is a $D$-space.
\endabstract

\endtopmatter

\document

\head 1. Introduction \endhead

 The class of $D$-spaces, introduced
by E. van Douwen in \cite{vDP}, is a very natural one.  $X$ is a
{\it $D$-space} iff, given a ``neighborhood assigment"
$\{N(x):x\in X\}$ (i.e, $x\in\text{Int} N(x)$ for each $x\in X$),
there is a closed discrete subset $D$ of $X$ such that
$X=\bigcup\{N(x):x\in D\}$.

There has been some interesting recent work on $D$-spaces due
especially to Arhangel'skii and Buzyakova \cite{AB},
Buzyakova\cite{Bu}, and Fleisner and Stanely \cite{FS}.  In
particular, Arhangel'skii and Buzyakova show that spaces having a
point-\newline countable base are $D$. Fleissner and Stanely
introduce the notion of $N$-sticky for a neighborhood assignment
$N$, a tool which simplifies many $D$-space arguments.  Buzyakova
obtained an interesting result in $C_p$-theory which illustrates
how $D$-spaces can be useful: she proved that $C_p(X)$ is
hereditarily $D$ for compact $X$. This can be viewed as an
``explanation" for the important, now classical, result of
Baturov\cite{Ba}, that Lindel\"of degree equals extent for
subspaces of these $C_p(X)$'s.

In the first part of this note, we introduce the notion of a
nearly good relation, and generalize the Fleissner-Stanley
$N$-sticky notion.  We observe that the point-countable base
result and the $C_p(X)$ result mentioned above follow easily from
general results about these notions.   Baturov's result holds more
generally for Lindel\"of $\Sigma$-spaces $X$, and Matveev asked if
$C_p(X)$ is hereditarily $D$ for such $X$. We exploit our general
results to obtain a positive answer to Matveev's question.

Another corollary of Buzyakova's result that $C_p(X)$ for compact
$X$ is hereditarily $D$ is that Eberlein compacta, which are
embeddable in such function spaces, are hereditarily $D$. This led
Arhangel'skii to ask if Corson compacts are hereditarily $D$. We
show that the answer is positive.  We  answer  another question of
Arhangel'skii on $D$-spaces by showing that a countably compact
space which is a finite union of $D$-spaces must be compact.
Finally, we solve a problem of van Douwen by showing that a
perfectly normal collectionwise-normal space of R. Pol is a
$D$-space.

The following notation will be used throughout: If $N$ is a
neighborhood assignment on $X$, and $D\subset X$, we let
$N(D)=\bigcup_{d\in D}N(d)$.

 \head 2. Nearly good sticky relations \endhead

 Let $X$ be a space.  We say that a relation $R$ on $X$ (resp.,
 from $X$ to $[X]^{<\o}$) is  {\it nearly
good} if $x\in \cl{A}$ implies $xRy$ for some $y\in A$ (resp.,
$xR\tilde{y}$ for some $\tilde{y}\in [A]^{<\o}$).

Further, if $N$ is a neighborhood assignment on $X$, $X'\subset
X$, and $D\subset X$, we say $D$ is {\it $N$-sticky mod R on $X'$}
if whenever $x\in X'$ and $xRy$ for some $y\in D$ (resp.,
$xR\tilde{y}$ for some $\tilde{y}\in [D]^{<\o}$), then $x\in
N(D)$. (In other words, it means that $N(D)$ contains all the
``relatives" of members (resp, finite subsets) of $D$ that are in
$X'$.)  We say more briefly that {\it $D$ is $N$-sticky mod $R$}
if $D$ is $N$-sticky mod $R$ on $X$.

For example, if $N$ is a neighborhood assignment and we define
$xRy\iff y\in N(x)$, then ``$N$-sticky mod $R$" means ``$x\in
N(D)$ whenever $N(x)\cap D\neq \0$ and is what Fleissner and
Stanley\cite{FS} called simply ``$N$-sticky" .  Obviously this $R$
is nearly good.

We begin with the following lemma, which is an immediate
consequence of the definitions.

\proclaim{Lemma 2.0} Let $N$ be a neighborhood assignment on $X$,
and $R$ a nearly good relation (on $X$, or from $X$ to
$[X]^{<\o}$).  If $D$ is $N$-sticky mod $R$ on $X'$, then
$\cl{D}\cap X'\subset N(D)$. \endproclaim

The next lemma will help us build closed discrete sets $D$ with
$X=N(D)$.

\proclaim{Lemma 2.1} \roster
 \item"(a)" Let $R$ be a nearly good relation on $X$.
 If for each $\a<\l$, $D_\a$ is a closed
discrete subset of $X\sm N(\bigcup_{\b<\a}D_\b)$ and $N$-sticky
mod $R$ on $X\sm N(\bigcup_{\b<\a}D_\b)$, then
$\bigcup_{\a<\l}D_\a$ is closed discrete. \item"(b)" Let $R$ be a
nearly good relation from $X$ to $[X]^{<\o}$. If for each $\a<\l$,
$D_\a$ is a closed discrete subset of $X\sm
N(\bigcup_{\b<\a}D_\b)$ and $\bigcup_{\b\leq\a}D_\b$ is $N$-sticky
mod $R$ on $X\sm N(\bigcup_{\b<\a}D_\b)$, then
$\bigcup_{\a<\l}D_\a$ is closed discrete.
 \endroster
\endproclaim

\demo{Proof} We prove (b) first.  Suppose $x$ is a limit point of
$\bigcup_{\a<\l}D_\a$. Since $R$ is nearly good, there are some
$\a'<\l$ and $\t{y}\in [\bigcup_{\b\leq \a'}D_{\b}]^{<\o}$ with
$xR\t{y}$. By the $N$-stickiness of $\bigcup_{\b\leq \a'}D_{\b}$ ,
we must have $x\in N(\bigcup_{\b\leq\a'}D_\b)$.  Let $\a$ be least
such that $x\in N(D_\a)$.  Then by the same argument, $x$ is not a
limit point of $\bigcup_{\b<\a}D_\b$. Since $N(D_\a)\cap
\bigcup_{\g>\a}D_\g=\0$, $x$ is not a limit point of $D$,
contradiction.

Part (a) is similar, noting that for relations on $X$ we only need
to apply $N$-stickiness to individual $D_\a$'s, instead of unions
of the type $\bigcup_{\b\leq \a}D_{\b}$.
 \qed
\enddemo

\proclaim{Proposition 2.2} Let $N$ be a neighborhood assignment on
$X$.  \roster\item"(a)" Suppose $R$ is a nearly good relation on
$X$ such that every non-empty closed subset $F$ of $X$ contains a
non-empty closed discrete    subset $D$ which is $N$-sticky mod
$R$ on $F$. Then there is a closed discrete $D^*$ in $X$ with
$N(D^*)=X$.
\item"(b)" Let $R$ be a nearly good relation from $X$ to
$[X]^{<\o}$.  Suppose that given any closed discrete $D$ and
non-empty closed $F\subset X\sm N(D)$ such that $D$ is $N$-sticky
mod $R$ on $F$, there is a non-empty closed discrete $E\subset F$
such that $D\cup E$ is $N$-sticky mod $R$ on $F$.  Then there is a
closed discrete $D^*$ in $X$ with $N(D^*)=X$.
\endroster

\endproclaim

\demo{Proof}  The proofs of (a) and (b) are essentially the same.
Inductively define closed discrete $D_\a\subset X\sm
N(\bigcup_{\b<\a}D_\b)$ satisfying the stickiness property given
by (a) or (b), until a stage $\l$ is reached such that
$X=N(\bigcup_{\a<\lambda}D_\a)$.  Then apply Lemma 2.1 to see that
$D^*=\bigcup_{\a<\l}D_\a$ is closed discrete.   \qed
\enddemo

By part (a) of this proposition, if we wish to prove that a
certain closed-hereditary property implies $D$, we just need to
prove that any neighborhood assignment $N$ on a space with the
property contains some non-empty $N$-sticky mod $R$ closed
discrete subset for some nearly good $R$ (as long as $R$ is
defined only in terms of $N$ and the property).

For example, suppose $X$ is left-separated and $N$ is a
neighborhood assignment on $X$.  W.l.o.g., $N(x)\subset
[x,\rightarrow)$, where the implied order is the order that
left-separates $X$.   Then every non-empty subset $F$ of $X$ has a
non-empty closed discrete $N$-sticky subset, namely the least
element of $F$.  So by Proposition 2.2(a), left-separated spaces
are $D$.  van Douwen and Pfeffer\cite{vDP} show this for the
so-called ``generalized left-separated" spaces, and this also
follows from Proposition 2.2(a) by a similarly easy argument.

In certain more complicated applications of Proposition 2.2(a), it
is natural to build a countable closed discrete $N$-sticky set.
Here the use of countable elementary submodels can significantly
simplify arguments.  At first glance, countable elementary
submodels do not seem as relevant to Proposition 2.2(b).  However,
they are relevant because it turns out that if 2.2(b) is true for
all countable $D$, it is true for all $D$.

\proclaim{Proposition 2.3} Let $N$ be a neighborhood assignment on
$X$ , and let $R$ be a nearly good relation from $X$ to
$[X]^{<\o}$. Suppose that given any countable closed discrete $D$
and non-empty closed $F\subset X\sm N(D)$ such that $D$ is
$N$-sticky mod $R$ on $F$, there is a countable non-empty closed
discrete $E\subset F$ such that $D\cup E$ is $N$-sticky mod $R$ on
$F$. Then there is a closed discrete $D^*$ in $X$ with $N(D^*)=X$.
\endproclaim

\demo{Proof} Suppose $D$ is closed discrete, $\0\neq F\subset X\sm
N(D)$ is closed, and $D$ is $N$-sticky mod $R$ on $F$.  We need to
show that there is a non-empty closed discrete $E\subset F$ such
that $D\cup E$ is $N$-sticky mod $R$ on $F$.   Then the required
$D^*$ exists by 2.2(b).

Let $D=\{d_\a<\a<\k\}$, where $\k=|D|$, and suppose any $D$ of
cardinality smaller than  $\k$ can be extended as required to an
$E$ such that $|E|\leq |D|+\o$.  By our assumption, $\k>\o$.
Inductively define non-empty closed discrete sets $E_\a\subset
F\sm N(\bigcup_{\b<\a}E_\b)$ of cardinality $\leq |\a|+\o$ such
that $\{d_\b:\b<\a\}\cup \bigcup_{\b\leq \a}E_\b$ is $N$-sticky
mod $R$ on $F$.  Stop the induction either at $\l=\k$, or at any
$\l<\k$ for which $F\sm N(\bigcup_{\a<\l}E_\a)=\0$.   It is easy
to see that $D\cup E$ is $N$-sticky mod $R$ on $F$, and, using
Lemma 2.1(b), that $E$ is closed discrete.  \qed
\enddemo

We now describe a general situation in which spaces can be shown
to be $D$-spaces because the hypothesis of Proposition 2.3 holds.

Given a neighborhood assignment $N$ on $X$, let us call a subset
$Z$ of $X$ {\it $N$-close} if $x,x'\in Z\imp x\in N(x')$
(equivalently, $Z\subset N(x)$ for every $x\in Z$).

\proclaim{Proposition 2.4} Let $N$ be a neighborhood assignment on
$X$.  Suppose there is a nearly good $R$ on $X$ (resp., from $X$
to $[X]^{<\o}$ ) such that for any $y\in X$ (resp., $\tilde{y}\in
[X]^{<\o}$), $R^{-1}(y)\sm N(y)$ (resp., $R^{-1}(\tilde{y})\sm
N(\tilde{y})$) is the countable union of $N$-close sets.  Then
there is a closed discrete $D$ such that $N(D)=X$. \endproclaim

{\bf Remark.}  Note that if $N$ and $R$ satisfy the hypotheses of
the proposition, then so does their restriction to any subspace.
So, if for any $N$ on $X$ we can produce such an $R$, then $X$ is
hereditarily $D$.

 \demo{Proof of the proposition}  We prove this in case $R$
 is a relation from $X$ to $[X]^{<\o}$.  By Proposition 2.3, we need
 only show that if $D$ is countable, closed discrete,
 and $N$-sticky mod $R$ on some non-empty closed
 $X'\subset X\sm N(D)$, then there is a non-empty countable
 closed discrete $E\subset X'$ such that $D\cup
 E$ is $N$-sticky mod $R$ on $X'$.

For $\tilde{y}\in [X]^{<\o}$, let $R^{-1}(\tilde{y})\sm
N(\tilde{y}) =\bigcup_{n\in\o}G_n(\tilde{y})$, where each
$G_n(\tilde{y})$ is $N$-close.  Put all relevant objects in a
countable elementary submodel $\Cal M$. Let $<_{\Cal M}$
well-order $\Cal M$ in type $\o$.  Choose $e_0\in X'\cap \Cal M$.
If $e_i\in X'\cap \Cal M$ has been defined for all $i<n$, look at
$$X_n'=\{x\in X'\sm N(\{e_i:i<n\}):xR\tilde{y}\text{ for some
}\tilde{y}\in [D\cup \{e_i:i<n\}]^{<\o}\}.$$

If $x\in X_n'$, then $x\in G_n(\tilde{y})\subset N(x)$ for some
$\tilde{y}\in [D\cup \{e_i:i<n\}]^{<\o}\}$.  Note any such
$G_n(\tilde{y})$ is in $\Cal M$.  Choose $e_n\in X_n'\cap \Cal M$
such that the corresponding $G_n(\tilde{y})$ is $<_{\Cal M}$ least
possible.

If $X_n'=\0$ for any $n>0$, then $D\cup \{e_i:i<n\}$ is closed
discrete and  $N$-sticky mod $R$ relative to $X'$ and we are done.
If $X_n'\neq \0$ for all $n>0$, let us show that if
$E=\{e_i:i<\o\}$, then $D\cup E$ is $N$-sticky mod $R$ on $X'$ and
closed discrete. Clearly $E$ is relatively discrete in $N(E)$, so
by Lemma 2.0,  it suffices to prove $D\cup E$ is $N$-sticky mod
$R$ on $X'$. To this end, suppose $x\in X'\sm N(D\cup E)$ and
$xR\tilde{y_0}$ for some $\tilde{y_0}\in [D\cup E]^{<\o}$. Then
for all sufficiently large $n$, we have $x\in X_n'$.  Let $n_0$ be
such that  $x\in G_{n_0}(\tilde{y_0})$, and note that
$G_{n_0}(\tilde{y_0})\in\Cal M$. Since $N(e_n)$ always contains
the $<_{\Cal M}$-least $G_n(\tilde{y})$ corresponding to some
$x\in X_n'$, eventually we chose  $e_n$ with $N(e_n)\supset
G_{n_0}(\tilde{y_0})$, which puts $x\in N(e_n)$, contradiction.
\qed
\enddemo

Recall that $X$ satisfies {\it open (G)} if for each $x\in X$ we
have a countable open neighborhood base  $\Cal B_x$ of $x$ such
that whenever $x\in \cl{A}$ and $N(x)$ is a neighborhood of $x$,
then for some $a\in A$ we have $x\in B\subset N(x)$ for some $B\in
\Cal B_a$.  Spaces with a point-countable base satisfy open (G),
but whether or not the reverse holds is an unsolved problem
\cite{CR}. We illustrate the use of Proposition 2.4 by proving the
following generalization of the Arhangel'skii-Buzyakova result
about point-countable bases (which can similarly be derived from
2.4).

 \proclaim{Proposition 2.5} Any space satisfying open (G) is a $D$-space.
\endproclaim
\demo{Proof} Let $X$ satisfy open (G), and let $N$ be a
neighborhood assignment.  Define
$$xRy\iff \exists B\in \Cal B_y\text{ with }x\in B\subset N(x).$$ It
is clear from the definition of open (G) that this $R$ is nearly
good.

 For each $B\in \Cal B_y$, let
$C(B)=\{x: x\in B\subset N(x)\}$. Then $C(B)$ is $N$-close, and
$R^{-1}(y)=\bigcup_{B\in \Cal B_y}C(B)$. By Proposition 2.4, $X$
is $D$. \qed
\enddemo

The framework encompassed by our Propositions 2.2(b),2.3, and 2.4
is implicit in Buzyakova's proof of the following, which we give
here as another illustration of the use of our Proposition 2.4.

 \proclaim{Proposition 2.6 \cite{Bu}} If $X$ is compact, then $C_p(X)$
 is hereditarily $D$. \endproclaim

\demo{Proof} Let $\Cal B$ be a countable base for the real line
$\Bbb R$.   For $S\subset X$ and $B\in \Cal B$, let $[S,B]=\{f\in
C(X):f(S)\subset B\}$.  For $A\subset C(X)$, let $\Cal G_A$ be the
set of all $G=\bigcap_{i<n}[S_i,B_i]$ where $B_i\in \Cal B$ and
$S_i$ can be written in the form $X\sm \bigcup_{a\in
A'}a^{-1}(B_a)$ for some finite $A'\subset A$.

Let $N$ be a neighborhood assignment on $C_p(X)$.  For $f\in
C_p(X)$ and $\tilde{g}\in [C_p(X)]^{<\o}$, define
$$fR\tilde{g}\iff \exists G\in \Cal G_{\tilde{g}}(f\in G\subset
N(f)).$$  Without using the terminology, Lemma 2.3 of \cite{Bu}
says exactly that this $R$ is nearly good. For each $G\in \Cal
G_{\tilde{g}}$, let $C(G)=\{f\in R^{-1}(\tilde{g}):f\in G\subset
N(f)\}$.  Note that $C(G)$ is $N$-close.  Since $\Cal
G_{\tilde{g}}$ is countable, we have that $R^{-1}(\tilde{g})$ is a
countable union of $N$-close sets.  So $C_p(X)$ is hereditarily
$D$ by Proposition 2.4. \qed
\enddemo

A similar use of Proposition 2.4 answers a question of Matveev
\cite{M}.

\proclaim{Proposition 2.7} Let $X$ be a Lindel\"of $\Sigma$-space.
Then $C_p(X)$ is hereditarily $D$. \endproclaim

\demo{Proof}  Since $X$ is Lindel\"of $\Sigma$, there are a cover
$\Cal K$ by compact sets and a countable collection $\Cal F$ such
that, whenever $K\in \Cal K$ and $K\subset U$, where $U$ is open,
then $K\subset F\subset U$ for some $F\in \Cal F$.

For $A\subset C(X)$, define $\Cal G_A$ just like in the proof of
Proposition 2.6, except that the $S_i$'s may have the form $F\sm
\bigcup_{a\in A'}a^{-1}(B_a)$ for some finite $A'\subset A$, where
$F\in \Cal F$.  Then define the relation $R$ just like before.
Since $\Cal F$ is countable, $\Cal G_A$ for countable $A$ is too,
so by the same argument each $R^{-1}(\tilde{g})$ is a countable
union of $N$-close sets.

Thus it remains to prove that $R$ is nearly good.  To this end,
suppose $f\in \cl{A}$ for some $A\subset C_p(X)$.  We need to show
that $f\in G\subset N(f)$ for some  $G\in \Cal G_A$.  (Note that
any $G\in \Cal G_A$ is in $\Cal G_{A'}$ for some finite $A'\subset
A$.) Since $\Cal G_A$ is closed under finite intersections, we may
assume $N(f)$ is a subbasic open set $[\{p\},B]$, where $p\in X$
and $B\in \Cal B$ .

Let $p\in K$, where $K\in \Cal K$.  Let $B'$ be open in $\Bbb R$
with $f(p)\subset B'\subset \cl{B'}\subset B$.  For each $y\in K$
with $f(y)\not\in B$, choose $B_y\in \Cal B$ containing $f(y)$
with $B_y\cap B'=\0$.  Since $f\in \cl{A}$, we can choose some
$a_y\in A$ with $a_y(y)\in B_y$ and $a_y(p)\in B'$.  By
compactness, there are $y_i$, $i<n$, such that the sets
$a_{y_i}^{-1}(B_{y_i})$ cover $K\sm f^{-1}(B)$.  Let $F\in \Cal F$
such that $K\subset F\subset f^{-1}(B)\cup
\bigcup_{i<n}a_{y_i}^{-1}(B_{y_i})$.  Let $S=F\sm \bigcup_{i<n}
a_{y_i}^{-1}(B_{y_i})$.  Then $[S,B]\in \Cal G_A$ and $f\in
[S,B]\subset [\{p\},B]=N(f)$.  \qed
\enddemo

\head 3. Corson compacts \endhead

A corollary of Bouziakova's result that $C_p(X)$ is hereditarily a
$D$-space whenever $X$ is compact is that Eberlein compact spaces
are hereditarily $D$. This prompted the natural question, due to
Arhangel'skii, whether Corson compact spaces are hereditarily $D$.
We will show that the answer is positive.

Recall that $X$ is Corson compact iff $X$ is compact and can be
embedded into a $\Sigma$-product of real lines.  Using the fact
that any closed interval in the real line containing $0$ is a
$\leq 2$-to-one continuous image of the Cantor set under a map $f$
with $f^{-1}(0)=0$, it is easy to see (and well-known) that any
Corson compact space is the continuous image of a $0$-dimensional
Corson compact space.   Also, the $D$-space property is preserved
by closed mappings \cite{BW}.  It follows that it suffices to
prove that $0$-dimensional Corson compact spaces are hereditarily
$D$.

\proclaim{Lemma 3.1}  Let $X$ be Corson compact and
$0$-dimensional. Then $X$ has a point-countable $T_0$-separating
cover $\Cal B$ by compact open sets which is closed under finite
intersections.
\endproclaim

\demo{Proof}  From $0$-dimensionality and the  fact that a compact
space $X$ is Corson compact iff $X$ has a point-countable
$T_0$-separating cover by open $F_\s$-sets \cite{MR}, it easily
follows that there is a point-countable $T_0$-separating
collection $\Cal B$ of compact open sets.  Take any such $\Cal B$
and close it under finite intersections. \qed
\enddemo

\proclaim{Lemma 3.2}  Let $\Cal B$  a point-countable
$T_0$-separating cover by compact open sets of a compact space $X$
which is closed under finite intersections. Then every $x\in X$
has a neighborhood base of sets of the form $B\sm \cup \Cal C$,
where $B\in \Cal B$ and  and $\Cal C$ is a  finite subcollection
of $\Cal B$.
\endproclaim

\demo{Proof}  Let $x\in U$, $U$ open.  For each $y\in X\sm U$,
either there is $B_y\in \Cal B$ with $x\in B_y$ and $y\not\in
B_y$, or there is $C_y\in \Cal B$ with $y\in C_y$ and $x\not\in
C_y$.  By compactness, there is a finite subcollection of
$\{C_y:y\in X\sm U\}\cup \{X\sm B_y:y\in X\sm U\}$ which covers
$X\sm U$.  Then take $B$ to be the intersections of the  $B_y$'s
from this finite subcover , and take $\Cal C$ to be the $C_y$'s.
\qed
\enddemo

Given a collection $\Cal S$ of finite sets, let $\Cal R(\Cal S)$
denote the collection of all roots of uncountable $\Delta$-systems
from $\Cal S$ (i.e., $R\in \Cal R(\Cal S)$ iff there is an
uncountable subcollection $\Cal S'$ of $\Cal S$ such that $S_0\cap
S_1=R$ whenever $S_0$ and $S_1$ are distinct elements of $\Cal
S'$). Then let  $$\Cal M(\Cal S)=\{R\in \Cal R(\Cal S):\not\exists
R'\in \Cal R(\Cal S)(R'\subsetneq R)\}\cup\{S\in \Cal
S:\not\exists R\in \Cal R(R\subseteq S)\}.$$

\proclaim{Lemma 3.3}  For any collection $\Cal S$ of finite sets,
the collection $\Cal M(\Cal S)$ is countable.
\endproclaim

\demo{Proof} Easy application of the $\Delta$-system lemma (that
any uncountable collection of finite sets contains an uncountable
$\Delta$-system).\qed
\enddemo

\proclaim{Theorem 3.4} Every Corson compact space is hereditarily
a $D$-space.  \endproclaim

\demo{Proof}  Let $X$ be Corson compact, and $Z\subset X$.  By the
remark preceding Lemma 3.1, we may assume $X$ is $0$-dimensional.
Then by Lemma 3.1, there is a point-countable $T_0$-separating
cover $\Cal B$ of $X$ consisting of compact open sets which is
closed under finite intersections.

Let $N(z)$, $z\in Z$, be a neighborhood assignment.   By Lemma
3.2, we may assume $N(z)=(B_z\sm \cup \Cal C_z)\cap Z$ for some
$B_z\in \Cal B$ and finite $\Cal C_z\subset \Cal B$. By our
observation after Proposition 2.2, we need only show that there
exists a non-empty closed discrete $N$-sticky subset $D$ of $Z$.
Recall that $N$-sticky means we need $D\cap N(z)\neq \0$ to imply
$z\in N(D)$.

To this end, put $X,\Cal B, Z, N,...$ in a countable elementary
submodel $M$ (of $H(\k)$ for some sufficiently large $\k$). Let
$\{\Cal B_i\}_{i\in\w}$ enumerate $ M\cap [\Cal B]^{<\w}$ in type
$\w$ such that each term is listed infinitely often.

At step $k$, $k\in \o$, we are going to define a finite subset
$F_k$ of $Z\cap M$ to put in $D$.  Look at $\Cal B_k$.  If $\Cal
B_k$ is not a singleton, let $F_k=\0$.  Otherwise, let $\Cal B_k=
\{B_k\}$, and  consider the collection
$$\Cal S_k=\{\Cal C:\exists z\in Z\cap B_k\sm N(\cup\{F_i:i<k\})
\text{ with }B_z=B_k\text{ and }\Cal C_z=\Cal C\}.$$  If $\Cal
S_k=\0$, let $F_k=\0$. Otherwise, continue as follows.  Note that
$\Cal S_k$ is in $M$ since all parameters in its definition are in
$M$. Thus $\Cal M(\Cal S_k)$ is in $M$, and since by Lemma 2.3 it
is countable, we have $\Cal M(\Cal S_k)\subset M$.

{\it Claim 1. For each $R\in \Cal M(\Cal S_k)$, there is a finite
set $F(R)\subset M\cap Z\cap B_k\sm N(\cup\{F_i:i<k\}) $ such that
$\bigcup_{z\in F(R)}N(z)=Z\cap B_k\sm \cup R$.}

To see this, consider $R\in\Cal M(\Cal S_k)$.  Let $$Z_k=\{z\in
Z\cap B_k\sm N(\cup\{F_i:i<k\}): B_z=B_k\}.$$  Note that if $R\in
\Cal S_k$, then $\Cal C_z=R$ for some $z\in Z_k$. By elementarity,
there is such a $z$, call it $z_R$, in $M$, and taking
$F(R)=\{z_R\}$ works. Suppose on the other hand that $R\not\in
\Cal S_k$.  Then $R$ is the root of some uncountable
$\Delta$-system  $\Cal S'\subset \Cal S_k$. Since $\Cal B$ is
point-countable, $\bigcap_{S\in \Cal S'}\cup(S\sm R)$ is empty. By
compactness, some finite subcollection of these $S\sm R$'s has
empty intersection. It follows that there is some finite subset
$F(R)$ of $Z_k$ such that $R\subset \Cal C_z$ for each $z\in F(R)$
and  $\bigcap_{z\in F(R)}\cup(\Cal C_z\sm R)=\0$.  Again by
elementarity, there is such an $F(R)$ in $M$, and it is clear that
this $F(R)$ satisfies the desired condition.

Having established Claim 1, let $R_k$ be the least member of $\Cal
M(\Cal S_k)$ in our indexing of $M\cap [\Cal B]^{<\o}$, and let
$F_k$ be the set $F(R)$ guaranteed by Claim 1 with $R=R_k$.

{\it Claim 2. If $j<k$ and $B_j=B_k$, then $R_j\neq R_k$.}

Suppose that $R_j=R_k$. By the construction, at stage $k$ there is
some $ z\in Z\cap B_k\sm N(\cup\{F_i:i<k\})$ with $R_k\subseteq
\Cal C_z $. So $z$ is in $B_k\sm \cup R_k=B_j\sm \cup R_j$ and is
not in $N(\bigcup_{i\leq j}F_i)$, contradicting that $F(R_j)$
satisfies the conclusion of Claim 1.

We let $D=\bigcup_{i\in\o}F_i$.  Note that by the construction, if
$z\in D$, then $N(z)\cap D$ is finite.  Thus $D$ is relatively
discrete.   It remains to prove that $D$ is closed and $N$-sticky,
 which follows easily from:

{\it Claim 3.  If $p\in Z$ and $B_p\in M$, then $p\in N(D)$.}

Suppose not.  Then for each $k\in \o$ such that $B_p=B_k$, the
following holds at stage $k$ of the inductive procedure for
building $D$:
$$\exists z\in Z\cap B_p\sm N(\bigcup_{i<k}F_i)\text{ with
}B_z=B_p\text{ and }\Cal C_z=\Cal C_p.$$

Thus $\Cal C_p\in \Cal S_k$, and so there is also some $R_k'\in
\Cal M(\Cal S_k)$ with $R_k'\subseteq \Cal C_p$.  Then there is
$R^*\subset \Cal C_p$ such that $R_k'=R^*$ for infinitely many
$k$.

 Recall that at a stage $k$ like this, the least $R\in \Cal
M(\Cal S_k)$ is selected and denoted by $R_k$.  Since there are
only finitely many possible $R$'s less than $R^*$, and $R^*$ has
the possibility of being selected infinitely often, it follows
from Claim 2 that $R^*=R_k$ for some $k$.  So for this $k$ we have
$F_k=F(R^*)$.   Then by Claim 1, $N(F_k)\supset Z\cap B_k\sm \cup
R^*$, which puts $p\in N(D)$, contradiction. \qed

\enddemo

\head 4. Finite unions of $D$-spaces and linearly Lindel\"ofness
\endhead

In the problems section of the Zoltan Balogh Memorial Topology
Conference booklet,  and also in \cite{A}, Arhangel'skii asked
whether the union of two $D$-spaces must be a $D$-space.  He also
asked if a countably compact space that is the union of two
$D$-spaces must be compact.
  In this section, we give a positive answer to the second
  question.  Our answer is a corollary to our more general result
  that any  space of countable extent which can be written as
   the finite union of $D$-spaces must be linearly Lindel\"of.

The first question, if it has a positive answer, would imply the
second (since countably compact $D$-spaces are compact), but that
one is still unsolved.  Another related problem from \cite{A} that
also remains unsolved is whether or not a countably compact space
that is a countable union of $D$-spaces must be compact.

Recall that a space $X$ is {\it linearly Lindel\"of} if every
increasing open cover of $X$ has a countable subcover.  This is
well-known to be equivalent to the statement that every subset of
$X$ of uncountable regular cardinality has a complete accumulation
point.  The following is another known characterization; for the
benefit of the reader, we include its easy proof.

\proclaim{Lemma 4.1} A space $X$ is linearly Lindel\"of iff
whenever $\Cal O$ is an open cover of $X$ of cardinality $\kappa$
and $\Cal O$ has no subcover of cardinality $<\k$, then
$cf(\k)\leq \w$.
\endproclaim

\demo{Proof} Suppose  $X$ is linearly Lindel\"of and  $\Cal
O=\{O_\a:\a<\k\}$ is an open cover of $X$ of cardinality $\kappa$
with no subcover of cardinality $<\k$.  Let $U_\a=\bigcup_{\b\leq
\a}O_\b$.  Then $\{U_\a:\a<\k\}$ is an increasing open cover,
which must therefore have a countable subcover
$\{U_{\a_n}\}_{n\in\o}$.  Since $\Cal O$ has no subcover of
cardinality $<\k$, $\{\a_n\}_{n\in\o}$ must be cofinal in $\k$.

For the other direction, suppose $X$ is not linearly Lindel\"of,
i.e., there is an  increasing open cover $\Cal U$  with no
countable subcover.  There is a cofinal subcollection $\Cal O$ of
$\Cal U$ of regular cardinality $\k$.  Note that $\Cal O$ has no
subcover of cardinality $<\k$ (by regularity of $\k$).   Since
$\Cal O$ has no countable subcover, we have $cf(\k)=\k>\w$. \qed
\enddemo

 \proclaim{Theorem 4.2} If $X$ has countable extent and can be
 written as the union of finitely many $D$-spaces, then $X$ is
 linearly Lindel\"of. \endproclaim

 \demo{Proof}.  Suppose $X$ satisfies the hypotheses, where
 $X=\bigcup_{i\leq k}X_i$ with each $X_i$ a $D$-space.  Suppose
 also by way of contradiction that $X$ is not linearly Lindel\"of
 and that $k$ is the least possible value for any counterexample
 to the theorem.  Of course $k>1$ since any $D$-space of
 countable extent is Lindel\"of.

By Lemma 4.1, there is an open cover $\Cal O=\{O_\a\}_{\a<\k}$ of
 some cardinality $\k$ with $cf(\k)>\w$ and such that $\Cal O$ has
 no subcover of cardinality $<\k$.  For each $x\in X$, let $\a_x$
 be least such that $x\in O_{\a_x}$ and consider the neighborhood
 assignment defined by $N(x)=O_{\a_x}$.

For each $i\leq k$, there is a relative closed discrete subset
$D_i$ of $X_i$ such that $\{N(d):d\in D_i\}$ covers $X_i$.  Since
$\Cal O$ has no subcover of smaller cardinality, there must be
some $i_0\leq k$ such that $|\{\a_d:d\in D_{i_0}\}|=\k$.  Note
that $Z=\cl{D_{i_0}}\sm D_{i_0}$ is closed in $X$ and is a subset
of $\bigcup_{i\neq i_0}X_i$. By minimality of $k$, $Z$ is linearly
Lindel\"of.  Applying this to the increasing open cover
$\{\cup_{\b<\a}O_\b:\a<\k\}$, there are $\a_n<\k$, $n\in\o$, such
that $\Cal U =\{\cup_{\b<\a_n}O_\b:n\in\o\}$ covers $Z$. Note that
$D_{i_0}\sm \cup \Cal U$ is closed discrete in $X$, so by
countable extent is countable.  By $cf(\k)>\w$, we have
$\delta=sup\{\a_n:n\in\o\}<\k$.  Hence there is some $d\in
D_{i_0}$ with $\a_d>\d$ and $d\in \cup \Cal U$.  But $d\in \cup
\Cal U$ implies $d\in O_\b$ for some $\b<\d$, whence $\a_d<\d$,
contradiction. \qed

\enddemo

\proclaim{Corollary 4.3}  Suppose $X$ is countably compact and a
finite union of $D$-spaces.  Then $X$ is compact. \endproclaim

\demo{Proof} Countably compact linearly Lindel\"of spaces are
compact.\qed
\enddemo

\head 5. Pol's space is  $D$ \endhead

In his talk at the International Conference in Topology in Matsue,
Japan, 2002, P.J. Nyikos mentioned the following problem related
to what he had called ``Classic Problem II" in the first volume
(1976) of Topology Proceedings: Is every (perfectly normal)
collectionwise-normal space with a point-countable base
paracompact?  This problem remains unsolved, not even consistency
results are known.   Arhangel'skii, recalling his result with R.
Buzyakova that spaces with a point-countable base are $D$-spaces,
asked in a verbal communication if it may even be that every
(perfectly normal) collectionwise-normal $D$-space is paracompact.
It turns out this essentially was asked earlier by van Douwen
\cite{vD}. He asked for a non-paracompact collectionwise-normal
space that is not ``trivially so".  He goes on to mention some
properties the space should have, and then says ``it would be even
better if the space is a $D$-space".  In this section we show that
that a perfectly normal collectionwise-normal non-paracompact
space constructed by R. Pol\cite{P} is a $D$-space, so this is an
example of the kind van Douwen asked for, and answers
Arhangel'skii's  question in the negative.

We use the following version of Pol's space $X$. For each $\a\in
\w_1$, choose a non-decreasing function $x_\a:\w\to \a$ with
$\a=sup\{x_\a(n):n\in \w\}$.  The set for $X$ is
$\{x_\a:\a<\o_1\}$.  For each $n\in \o$ and $\s\in \o_1^{n}$, let
$[\s]=\{x\in X:x\r n=\s\}$.  Then for each $\a<\o_1$ and $n\in
\o$, let $B(\a,n)=\{x_\b:\b\leq \a\text{ and }x_\b\r n=x_\a\r
n\}$. Note that  $B(\a,n)=[x_\a\r n]\cap \{x_\b:\b\leq \a\}$. The
$B(\a,n)$'s form a basis for Pol's topology on $X$, which is
clearly finer than the metric topology generated by the $[\s]$'s,
and is also finer than the ``interval" topology generated by sets
of the form $\{x_\g:\a<\g\leq \b\}$, where $\a,\b\in\o_1$.

\proclaim{Theorem 5.1} Pol's space $X$ described above is a
$D$-space.
\endproclaim

\demo{Proof}  Recall that a non-stationary subset $A$ of $\o_1$ is
metrizable; similarly, $X_A=\{x_\a:\a\in A\}$ is metrizable
whenever $A$ is non-stationary.

Another fact about $X$ we shall use is that every uncountable
subset of $X$ contains an uncountable closed discrete set.  To see
this, note that since the topology is finer than the interval
topology, every uncountable subset has an uncountable relatively
discrete subset; then apply perfect normality.

Now suppose we are given an open neighborhood assignment for $X$.
W.l.o.g., this can be coded by $f:\w_1\to \w$, where $B(\a,f(\a))$
is the assigned open neighborhood of $x_\a$.

Let $\Sigma$ denote all $\s\in \o_1^{<\w}$ satisfying: \roster
\item"(i)"$\exists$ stationary $S_\s$ such that $x_\a\r f(\a)=\s$
for all $\a\in S_\s$;
\item"(ii)" $\s$ is minimal w.r.t. (i) (i.e., no proper initial
segment of $\s$ satisfies (i)).
\endroster

Let $A=\{\b\in \o_1:x_\b\not\in \bigcup_{\s\in \Sigma}[\s]\}$, and
let $X_A=\{x_\b:\b\in A\}$.  Then $X_A$ is closed in $X$.  Also,
an easy pressing-down argument shows that $A$ is non-stationary
and hence $X_A$ is metrizable.   Thus there is a closed discrete
subset $D_0$ of $X_A$ such that $\{B(\b,f(\b)):x_\b\in D_0\}$
covers $X_A$.

Let $U=\bigcup \{B(\b,f(\b)):x_\b\in D_0\}$.

{\it Claim 1. If $\s\in \Sigma$ and $[\s]\not\subset U$, then for
sufficiently large $\a\in S_\s$, $x_\a\not\in U$.}

To prove  Claim 1, suppose that $x_\a\in U$ for unboundedly many
$\a\in S_\s$.   Let $x_\g\in [\s]$.  Consider $\a\in S_\s$ with
$\a>\g$ and $x_\a\in U$.  Then $x_\a\in B(\b,f(\b))$ for some
$x_\b\in D_0$. Note that $\a\leq \b$.  Since $x_\b\not\in
[\s]=[x_\a\r f(\a)]$ and $x_\a\r f(\b)=x_\b\r f(\b)$, it must be
the case that $\a<\b$ and $f(\a)>f(\b)$.  Then $x_\g\in [\s]\cap
\{x_\d:\d\leq \a\}\subset [x_\b\r f(\b)]\cap \{x_\d:\d\leq
\b\}=B(\b,f(\b))\subset U$.  Hence $[\s]\subset U$, which proves
Claim 1.

For each $\s\in \Sigma$ with $[\s]\not\subset U$, by Claim 1 and
the fact that every uncountable subset of $X$ contains an
uncountable closed discrete set, there exists an unbounded
$T_\s\subset S_\s$ such that $x_\a\not\in U$ for any $\a\in T_\s$
and $E_\s=\{x_\a:\a\in T_\s\}$ is closed discrete.  Let
$D=D_0\cup\bigcup\{E_\s:\s\in \Sigma, [\s]\not\subset U\}$.

{\it Claim 2.  $D$ is closed discrete.} Let $x\in X$.  If $x\in
U$, then $x$ is not a limit point of $D$ since $U$ misses all the
$E_\s$'s, and $D_0$ is closed discrete.  On the other hand, if
$x\not\in U$, then there is a unique $\s\in \Sigma$ with $x\in
[\s]$, and $[\s]$ misses $D_0$ and all $E_\tau$'s  with $\tau\in
\Sigma$ and $\tau\neq\s$.

The next claim completes the proof of the example.

 {\it Claim 3. $\{B(\a,f(\a)):x_\a\in D\}$  covers $X$.}
Let $x\in X$.  If $x\in U$, we are done, so suppose $x\not\in U$.
Then $x\not\in X_A$, so there is a unique $\s\in \Sigma$ with
$x\in [\s]$.  Then $[\s]\not\subset U$, so $T_\s$ and $D_\s$ are
defined.  Say $x=x_\g$.  Choose $\a\in T_\s$ with $\a>\g$. Then
$x_\a\in D$ and $x=x_\g\in [\s]\cap \{x_\b:\b\leq
\a\}=B(\a,f(\a))$, which completes the proof. \qed
\enddemo

\widestnumber\key{BGTW}

\tenpoint 


\Refs

\tenpoint

\ref  \key{A} \by A.V. Arhangel'skii \paper $D$-spaces and finite
unions \jour preprint \endref

\ref\key{AB} \by A.V. Arhangel'skii and R. Buzyakova \paper
Addition theorems and $D$-spaces \jour Comment. Math. Univ. Car.
\vol 43 \yr 2002 \pages 653--663 \endref

\ref\key{Ba} \by D. Baturov \paper On subspaces of function spaces
\jour Vestnik MGU Mat. Mech. \vol 4 \yr 1987 \pages 66-69 \endref

\ref\key{BW} \by C.R. Borges and A.C. Wehrly \paper A study of
$D$-spaces \jour Topology Proc. \vol 16 \yr 1991 \pages 7-15
\endref


\ref\key{Bu} \by R. Buzyakova \paper Hereditary $D$-property of
function spaces over compacta \jour preprint \endref

\ref\key{CR} \by P.J. Collins and G.M. Reed \paper The
point-countable base problem \inbook Open Problems in Topology
\eds J. van Mill and G.M. Reed \publ Elsevier (North-Holland)
\publaddr Amsterdam \yr 1990 \pages 237-259 \endref


\ref\key{vD} \by E.K. van Douwen \paper Problem D-34 in the
Problems Section \jour Topology Proc. \vol 7 \yr 1982 \pages 383
\endref

 \ref\key{vDP} \by E.K. van Douwen and W. Pfeffer \paper Some
properties of the Sorgenfrey line and related spaces \jour Pacific
J. Math. \vol 81 \yr 1979 \pages 371-377
\endref


 \ref\key{FS} \by W. Fleissner and A. Stanley \paper
$D$-spaces \jour Topology Appl.\vol 114 \yr 2001 \pages 261--271
\endref

\ref\key{M} \by M. Matveev \paper $C_p(X)$ and $D$-spaces \jour
Zoltan Balogh Memorial Topology Conference Contributed Problems
\vol http://notch.mathstat.muohio.edu/ \newline
balogh\_conference/all\_prob.pdf \yr 2000  \pages 11 \endref


\ref\key{MR} \by E. Michael and M.E. Rudin \paper A note on
Eberlein compacts \jour  Pacific J. Math. \vol 72 \yr 1977 \pages
487--495\endref

\ref\key{P} \by R. Pol \paper A perfectly normal locally
metrizable non-paracompact space\jour Fundamenta Math. \vol 97 \yr
1977 \pages  37--42 \endref

\endRefs
\enddocument